# Combinatorial Ricci Flows on Surfaces

## Bennett Chow and Feng Luo

### Abstract

We show that the analog of Hamilton's Ricci flow in the combinatorial setting produces solutions which converge exponentially fast to Thurston's circle packing on surfaces. As a consequence, a new proof of Thurston's existence of circle packing theorem is obtained. As another consequence, Ricci flow suggests a new algorithm to find circle packings.

### §1. Introduction

1.1. For a compact surface with a Riemannian metric $(X, g_{ij})$, R. Hamilton in [Ha] introduced the 2-dimensional Ricci flow defined by the equation $\frac{dg_{ij}}{dt} = -2Kg_{ij}$ where $K$ is the Gaussian curvature of the surface. It is proved in [Ha] and [Cho] that for any closed surface with any initial Riemannian metric, the solution of the Ricci flow exists for all time, and after normalizing the solution to have a fixed area, the solution converges to the constant curvature metric conformal to the initial metric as time goes to infinity. The purpose of the paper is to study the analogous flow in the combinatorial setting. The simplest case of our setup is as follows. Given a triangulated surface, we assign to each vertex a positive radius and realize the surface by a piecewise flat cone metric with singularities at the vertices where the length of an edge is the sum of the two radii associated to the end points. We call these the (tangential type) circle-packing metrics based on the triangulation. We consider a natural flow (the combinatorial Ricci flow) on the space of all tangential type circle-packing metrics. Our results give a complete description to the asymptotic behavior of the solution of the Ricci flow. If we use hyperbolic geometry as the background metric, then for compact surfaces of negative Euler characteristic, we show that the combinatorial Ricci flow has solutions for all time for any initial metric and converges exponentially fast to the circle packing metric constructed by Thurston.

1.2. For simplicity, we shall state our results for triangulations on compact surfaces without boundary. For the most general version involving cellular-decompositions, see §2, §5 and §6.

Suppose $X$ is a closed surface and $T$ is a triangulation on $X$. Let $V = \{v_1, v_2, ..., v_N\}$ be the set of vertices in $T$. The set of all edges and triangles in $T$ are denoted by $E$ and $F$. Unless mentioned otherwise a *weight* on the triangulation is a function $\Phi: E \to [0, \pi/2]$. (In lemma 2.3, proposition 3.1 and proposition 4.1, we allow the weight to be $\Phi: E \to [0, \pi)$).

Given a weighted triangulation $(T, \Phi)$ of a compact surface $X$, we consider the circle packing metrics on the surface $X$ based on $(T, \Phi)$ defined as follows.

To each vertex $v_i$ in the triangulation $T$, assign a positive number $r_i$. Now realize each edge $e_{ij}$ joining $v_i$ to $v_j$ by a Euclidean segment of length

$$l_{ij} = \sqrt{r_i^2 + r_j^2 + 2r_ir_j \cos(\Phi(e_{ij}))}.$$

Note that, due to $\Phi(e) \leq \pi/2$, these positive numbers $\{l_{ij}, l_{jk}, l_{ki}\}$ satisfy triangle inequalities (see [Th, lemma 13.7.2]). (Also note that if the weight $\Phi = 0$, then $l_{ij} = r_i + r_j$. This corresponds to the case that circle patterns come from tangency.) Thus we can realize each triangle $\Delta v_i v_j v_k$ by a Euclidean triangle of edge lengths $l_{ij}, l_{jk}, l_{ki}$. The triangle is formed by the centers of three circles of radii $r_i$, $r_j$ and $r_k$ intersecting at angles $\Phi(e_{ij}), \Phi(e_{jk})$ and $\Phi(e_{ki})$ where $e_{nm}$ is the edge joining $v_n$ and $v_m$. This produces a flat cone metric on the surface $X$ with singularities at the vertices. We call these metrics the circle-packing metrics (based on $(T, \Phi)$). Let $a_i$ be the cone angle at the vertex $v_i$ which is the sum of all inner angles having vertex $v_i$. The *curvature* $K_i$ at $v_i$ is defined to be $2\pi - a_i$. The *combinatorial Ricci flow* is the following,

$$(1.1) \qquad \frac{dr_i}{dt} = -K_i r_i$$

It is easy to show that for any initial circle-packing metric, the solution of (1.1) with the given initial value exists for all time $t \geq 0$. Our goal is to understand the asymptotic behavior of the solution. To this end, it is useful to consider the normalized Ricci flow,

$$(1.2) \qquad \frac{dr_i}{dt} = -(K_i - K_{av}) r_i$$

where $K_{av} = 1/N \sum_{i=1}^{N} K_i$ is the average curvature. By the Gauss-Bonnet theorem (proposition 3.1), $K_{av} = 2\pi \chi(X)/N$ is time independent. The two equations (1.1) and (1.2) are related by a change of scales. Namely $r_i(t)$ is a solution to (1.1) if and only if $e^{2\pi t \chi(X)/N} r_i(t)$ is a solution of (1.2). Furthermore, if $r_i(t)$ is a solution to the normalized equation (1.2), the product $\prod_{i=1}^{N} r_i(t)$ is a constant. For this reason, when we consider the asymptotic behavior of the Ricci flow, it is more convenient to deal with the normalized Ricci flow (1.2). Let $\mathbf{R}_{>0}$ be the set of all positive real numbers.

**Definition.** A solution to (1.2) is called *convergent* if
 (1) $\lim_{t \to \infty} K_i(t) = K_i(\infty) \in (-\infty, 2\pi)$ exists for all $i$, and
 (2) $\lim_{t \to \infty} r_i(t) = r_i(\infty) \in \mathbf{R}_{>0}$ exists for all $i$.
 A convergent solution is called *convergent exponentially fast* if there are positive constants $c_1, c_2$ so that for all time $t \geq 0$,

$$|K_i(t) - K_i(\infty)| \leq c_1 e^{-c_2 t},$$

and

$$|r_i(t) - r_i(\infty)| \leq c_1 e^{-c_2 t}.$$

Given any subset $I \subset V$, let $F_I$ be the set of all cells in $T$ whose vertices are in $I$ and let the link of $I$, denoted by $Lk(I)$, be the set of pairs $(e, v)$ of an edge $e$ and a vertex $v$ so that (1) the end points of $e$ are not in $I$ and (2) the vertex $v$ is in $I$ and (3) $e$ and $v$ form a triangle.

**Theorem 1.1.** *Suppose $(T, \Phi)$ is a weighted triangulation of a closed connected surface $X$. Given any initial circle-packing metric based on the weighted triangulation, the solution to the normalized combinatorial Ricci flow (1.2) in the Euclidean geometry with the given initial value exists for all time and converges if and only if for any proper subset $I \subset V$,*

$$(1.3) \qquad 2\pi |I| \chi(X)/|V| > - \sum_{(e,v) \in Lk(I)} (\pi - \Phi(e)) + 2\pi \chi(F_I).$$

*Furthermore, if the solution converges, then it converges exponentially fast to the metric of constant curvature $2\pi \chi(X)/N$.*

By the work of Thurston on circle packing, condition (1.3) is equivalent to the existence of constant curvature circle packing metric based on the weighted triangulation on the surface.

In the case that the surface has negative Euler characteristic, it is natural to replace Euclidean triangles by the hyperbolic ones. The combinatorial Ricci flow in the hyperbolic geometry is the following. Given a surface $X$ with a weighted triangulation $(T, \Phi)$ where $\Phi(E) \subset [0, \pi/2]$, we assign to each vertex of T a positive radius $r_i$. We realize each triangle $\Delta v_i v_j v_k$ by a hyperbolic triangle formed by the centers of three circles of radii $r_i$, $r_j$ and $r_k$ intersecting at angles $\Phi(e_{ij})$, $\Phi(e_{jk})$ and $\Phi(e_{ki})$. Let $a_i$ be the cone angle at $v_i$ and $K_i = 2\pi - a_i$ be its curvature. The combinatorial Ricci flow in the hyperbolic metric is,

$$(1.4) \qquad \frac{dr_i}{dt} = -K_i \sinh r_i.$$

**Theorem 1.2.** *Suppose $(T, \Phi)$ is a weighted generalized triangulation of a closed connected surface $X$ of negative Euler characteristic. Given any initial choice of circle-packing metric, the solution to (1.4) in the hyperbolic background geometry with the given initial value exists for all time and converges if and only if the following two conditions hold.*

*(1.5a) for any three edges $e_1, e_2, e_3$ forming a null homotopic loop in $X$, if $\sum_{i=1}^{3} \Phi(e_i) \geq \pi$, then these three edges form the boundary of a triangle of $T$, and*

*(1.5b) for any four edges $e_1, e_2, e_3, e_4$ forming a null homotopic loop in $X$, if $\sum_{i=1}^{4} \Phi(e_i) \geq 2\pi$, then these $e_i$'s form the boundary of the union of two adjacent triangles.*

*Furthermore, if it converges, then it converges exponentially fast to a hyperbolic metric on $X$ so that all vertex angles are $2\pi$.*

Conditions (1.5a) and (1.5b) were obtained by Thurston in his work on circle packing on surfaces. Indeed, according to [Th] theorem 13.7.1, these two conditions are equivalent to the existence of circle packing metric with curvature equal to zero at each vertex.

The following result shows that condition (1.3) is almost the same as conditions (1.5a) and (1.5b) for surfaces with non-negative Euler characteristic.

**Proposition 1.3** *(a) If the Euler characteristic of the surface $X$ is positive, then condition (1.3) is valid for all triangulations of $X$ so that there are no null homotopic loops consisting of 3 edges $e_1, e_2, e_3$ with $\sum_{i=1}^{3} \Phi(e_i) \geq \pi$.*

*(b) If the Euler characteristic of the surface $X$ is zero, then (1.3) is valid for all triangulations of $X$ so that (1.5a) and (1.5b) hold.*

Note that if the weight is the zero function, the conditions in theorem 1.2 and proposition 1.3 automatically hold.

1.3. The exponential convergence of the Ricci flow may have some applications in practice. Recall that in his original proof of the existence of zero curvature circle packing metric, Thurston gave an algorithm to find the circle packing. This algorithm adjusts the radii at vertices one at a time. It is shown in [Mo] that Thurston's algorithm converges in polynomial time to the circle packing. To the best of our knowledge, the existing computer softwares (for instance, [St]) are based on Thurston's algorithm. It seems Ricci flow is more natural and may produce a faster algorithm to find circle packing metrics.

1.4. In the case that the background metric is spherical, we consider triangulations with zero weights. The combinatorial Ricci flow is

$$(1.5) \qquad \frac{dr_i}{dt} = -K_i \sin r_i$$

where we realize each triangle $\Delta v_i v_j v_k$ by a spherical triangle of edge lengths $r_i + r_j, r_j + r_k, r_k + r_i$ (note that we must impose $r_i \in (0, \pi)$ and $r_i + r_j + r_k < \pi$ in this case). Numerical simulations by L.T. Cheng [Che] shows that the situation is very complicated even in the simplest case of triangulations of the 2-sphere, i.e., the boundary of a tetrahedron. For instance, if we start with the boundary of the regular spherical tetrahedron, then the metric will shrink to a point in finite time. If we rescale it to have a fixed area, the cone angles will tend to $\pi$ (i.e., tends to the boundary of the Euclidean regular tetrahedron if rescaled). But in some other case, there are initial metrics so that as time goes to a finite number, the metric tends to the round 2-sphere. We do not know what the asymptotic behavior of the solution should be in general.

1.5. The proofs of the above theorems follow the same strategy used by Hamilton in his original approach in the smooth case. Namely, one first derives the evolution equation for the curvature. These equations are essentially of parabolic type. Thus the maximum principle can be applied to control the behaviors of the curvature. In some senses, what we are doing are quantifications of Thurston's original proof. For instance, the fact that the curvature evolution is a heat equation is essentially a consequence of Thurston's lemma 13.7.3 in [Th]. We also observed the fact that $\partial \theta_i / \partial r_j \sinh(r_j) = \partial \theta_j / \partial r_i \sinh(r_i)$ (for hyperbolic case) and $\partial \theta_i / \partial r_j r_j = \partial \theta_j / \partial r_i r_i$ (for Euclidean case) where $\theta_i$ is the inner angle at i-th vertex corresponding to $r_i$ radius circle (lemma 2.3). This was first observed in the paper of Colin de Verdiere [CV] for the zero weighted triangulation. This makes the Ricci flow variational. Namely, it is essentially a gradient flow of a convex function.

We observed that the degeneration condition obtained in Thurston's work implies that Palais-Smale condition holds for the convex functional. Thus establishes the main results. We remark that the convexity was first observed in [CV] in the case of the zero weighted triangulation. Also if a triangulation has zero weight, i.e., $\Phi = 0$, then the convergence part of theorem 1.2 can also be deduced from the results of Thurston [Th] and Colin de Verdiere [CV]. Our proofs are different from that of Colin de Verdiere [CV]. Colin de Verdiere used a result in network flow to produce the coherent system of angles. From this system, he then produced the convex function and showed that it is proper. In our case, it is the maximum principle and Palais-Smale condition which imply the exponential convergence. Another related work is [BS].

1.6. The paper is organized as follows. In §2, we collect some known results concerning triangles. Most of the proofs will be deferred to the appendix. In §3, we derive the evolution equation for the curvature and discuss some of its consequences. In §4, we establish a compactness theorem concerning the solutions of the Ricci flow. In §5, we prove theorem 1.2. In §6, we prove theorem 1.1. In the appendix, we derive some basic properties of the triangles.

1.7. The work is supported in part by the NSF DMS 0203926 (at UCSD) and DMS 0103843 (at Rutgers University).

## §2. Preliminaries on Triangles

Let $K^2$ be one of the three geometries, $E^2$ (Euclidean), $S^2$ (spherical) and $H^2$ (hyperbolic), i.e., a simply connected surface with a complete Riemannian metric of curvature equal to 0, 1, and -1. For simplicity, we define a function $s(x)$ according the geometry $K^2$ as follows. If $K^2 = E^2$, then $s(x) = x$; if $K^2 = H^2$, then $s(x) = \sinh(x)$; and if $K^2 = S^2$, then $s(x) = \sin(x)$.

A function is called *elementary* if it can be obtained from the polynomials, the exponential functions, logarithmic functions, the trigonometric functions and their inverses by the five operations of addition, subtraction, multiplication, division and composition.

2.1. The following lemma of Thurston guarantees the existence of a circle pattern in the geometry $K^2$.

**Lemma 2.1.** ([Th], lemma 13.7.2) *For any three angles $\Phi_i, \Phi_j, \Phi_k \in [0, \pi/2]$ and any three radii $r_i, r_j, r_k \in \mathbf{R}_{>0}$, there exists a configuration of three circles in both Euclidean and hyperbolic geometries, unique up to isometry, having radii $r_i, r_j, r_k$ and meeting in angles $\Phi_i, \Phi_j, \Phi_k$. Furthermore, if $r_i + r_j + r_k < \pi$, then such a configuration of three circles also exists in the spherical geometry $S^2$ and is unique up to isometry.*

Indeed, let $x_i$ be the distance between the centers of the $j$-th and $k$-th circles. Then it suffices to verify that $\{x_i, x_j, x_k\}$ satisfies the triangle inequalities that sum of two is larger than the third and also $x_i + x_j + x_k < 2\pi$ in the case of $S^2$. To see the first, we have $x_i \leq r_j + r_k \leq x_j + x_k$ since $\Phi_n \leq \pi/2$. Furthermore, the equality $x_i = x_j + x_k$ does not hold. The second inequality follows from $x_i + x_j + x_k \leq 2(r_i + r_j + r_k)$.

Consider the triangle formed by the centers of the three circles in lemma 2.1 in $K^2$. Let $\theta_i$ be the inner angle at the i-th vertex corresponding to the radius $r_i$ circle. Then $\theta_i$ is a smooth function of the radii $r_i, r_j, r_k$. The following is known in Thurston's notes.

**Lemma 2.2.** ([Th], lemma 13.7.3). *In all three geometries $E^2, H^2$ or $S^2$, if all angles $\Phi_n$'s are in $[0, \pi/2]$, then we have*
   (a) $\partial \theta_i / \partial r_i < 0$,
   (b) $\partial \theta_i / \partial r_j > 0$ for $i \neq j$,
   (c) $\partial(\theta_i + \theta_j + \theta_k)/\partial r_i < 0$ in $H^2$, $\partial(\theta_i + \theta_j + \theta_k)/\partial r_i > 0$ in $S^2$ and $\partial(\theta_i + \theta_j + \theta_k)/\partial r_i = 0$ in $E^2$.
*Furthermore, the partial derivatives $\partial \theta_n / \partial r_m$ are elementary functions in $r_i, r_j$ and $r_k$ where $n, m \in \{i, j, k\}$.*

Note part (c) follows from Thurston's geometric argument that the area of the triangle increases when $r_i$ increases. A computational proof of lemma 2.2 which gives the explicit formulas and the geometric meaning of these partial derivatives can be obtained using lemma A1 in the appendix. See the proof of lemma 2.3.

The next lemma was first obtained by Colin de Verdiere [CV] in the case where all angles $\Phi_n$'s are zero.

**Lemma 2.3.** *In all three geometries $E^2, H^2$ and $S^2$ and for all angles $\Phi_i, \Phi_j, \Phi_k \in [0, \pi)$ so that the triangle of radii $r_i$, $r_j$ and $r_k$ intersecting at these angles exists,*

$$\partial \theta_i / \partial r_j s(r_j) = \partial \theta_j / \partial r_i s(r_i).$$

We defer the proof to the appendix. In the appendix, we give two proofs. The simpler proof uses the mathematical software Maple VII.

One consequence of the lemmas is the following.

**Proposition 2.4.** *Suppose all angles $\Phi_n$ are in $[0, \pi/2]$. If the radii $r_n$'s evolve in time according to some differential equation $dr_n(t)/dt = -K_n s(r_n)$ where $K_n$ is a function of $r_i, r_j$ and $r_k$ and $n = i, j, k$, then the evolution of the inner angle $\theta_i$ satisfies the equation,*

$$d\theta_i/dt = -A_{ij}(K_j - K_i) - A_{ik}(K_k - K_i) - A_i K_i \lambda$$

*where $A_{rs}$ for $r \neq s$ and $A_n$ are positive valued elementary functions in $r_i, r_j$ and $r_k$, $\lambda$ is the curvature of the geometry $K^2$. Furthermore $A_{rs} = A_{sr}$.*

**Proof.** By the chain rule, we have

$$d\theta_i/dt = \partial \theta_i / \partial r_i r'_i + \partial \theta_i / \partial_j r'_j + \partial \theta_i / \partial r_k r'_k$$

$$= -\partial\theta_i/\partial r_i s(r_i) K_i - \partial\theta_i/\partial r_j s(r_j) K_j - \partial\theta_i/\partial r_k s(r_k) K_k$$

$$= -\partial\theta_i/\partial r_j s(r_j)(K_j - K_i) - \partial\theta_i/\partial r_k s(r_k)(K_k - K_i) - D_i K_i.$$

Here $D_i = \partial\theta_i/\partial r_i s(r_i) + \partial\theta_i/\partial r_j s(r_j) + \partial\theta_i/\partial r_k s(r_k)$. By lemma 2.3, we can write $D_i$ as $\sum_{l \in \{i,j,k\}} \partial\theta_l/\partial r_i s(r_i) = s(r_i)\partial(\sum_{l \in \{i,j,k\}} \theta_l)/\partial r_i$. By lemma 2.2 (c), the last partial derivative is negative in the case $K^2 = H^2$, is positive in the case $K^2 = S^2$ and is zero if $K^2 = E^2$. Define $A_i = \lambda D_i$ when $K^2 = S^2$ or $H^2$ and $A_i = 1$ when $K^2 = E^2$. Then $A_i > 0$. Let $A_{mn} = \partial\theta_m/\partial r_n s(r_n)$ for $n \neq m$. Then by lemma 2.2, it is positive and is an elementary function. Furthermore, the fact $A_{mn} = A_{nm}$ follows from lemma 2.3. Thus the result follows.

## §3. The 2-Dimensional Combinatorial Ricci Flow

3.1. *Basic setup and the Gauss-Bonnet formula*

We set up the Ricci flow in the most general setting as follows.

Suppose $X$ is a connected compact surface with possibly non-empty boundary. A *generalized triangulation $T$* of $X$ is a cell-division of $X$ into cells so that the lifting of the cell-division in the universal cover is a triangulation. In particular, each 2-cell (triangle) lifts to an embedding in the universal cover and there are no vertices of degree 1 or 2. Also there are no null homotopic loops formed by at most 2 edges. Furthermore, given three vertices, there is at most one triangle having them as vertices. The set of all vertices (0-cells), edges (1-cells) and triangles (2-cells) in $T$ are denoted by $V, E, F$. Unless mentioned otherwise, a *weight* on the generalized triangulation $T$ is a function $\Phi : E \to [0, \pi/2]$. The pair $(T, \Phi)$ is called a *weighted* generalized triangulation.

This definition of generalized triangulation is slightly different from Theorem 13.7.1 in Thurston's note. In fact, we were not able to rule out the degree 2 vertices from Thurston's definition. On the other hand, it is clear that if a triangulation is realized by a circle packing, there cannot be any degree two vertices.

Fix a background geometry $K^2$. Given a compact surface $X$ together with a weighted generalized triangulation $(T, \Phi)$, we assign to each vertex $v_i \in V = \{v_1, ..., v_N\}$ a positive number $r_i$. If the geometry $K^2$ is spherical $S^2$, we assume that $r_i + r_j + r_k < \pi$ holds for all triples $\{r_i, r_j, r_k\}$ so that the corresponding vertices form a triangle in $T$.

For each triangle $\Delta v_i v_j v_k$ in $F$ having edges $e_{ij}, e_{jk}$ and $e_{ki}$, using lemma 2.1 we can realize the triangle as the one formed by the centers of the three circles of radii $r_i, r_j, r_k$ meeting in angles $\Phi(e_{ij})$, $\Phi(e_{jk})$ and $\Phi(e_{ki})$ in geometry $K^2$. This produces a Riemannian metric on the surface $X$ of constant curvature having singularities at the vertices. Let $a_i$ be the sum of all inner angles having vertex $v_i$ in triangles. Then the *curvature* at $v_i$ is defined to be $K_i = 2\pi - a_i$ when the vertex $v_i$ is in the interior of $X$ and $K_i = \pi - a_i$ when $v_i$ is in the boundary.

Our first observation is the following,

**Proposition 3.1.** (Gauss-Bonnet) *Let $\lambda = -1, 0, 1$ be the curvature of one of the three geometries $H^2, E^2, S^2$ used in the construction. If the circle packing metric based on a* *weighted generalized triangulation* $(T, \Phi)$ *exists where* $\Phi : E \to [0, \pi)$, *then the total curvature*

$$\sum_{i=1}^{N} K_i = 2\pi\chi(X) - \lambda Area(X).$$

Note that the three terms in the Gauss-Bonnet formula are additive with respect to gluing along surface boundary. Thus by doubling the surface $X$ across the boundary if $X$ has non-empty boundary, it suffices to show the result for compact surfaces $X$ without boundary. If $X$ has no boundary, let $\theta_i^{jk}$ be the inner angle of the triangle $\Delta v_i v_j v_k$ at $v_i$ ($\theta_i^{jk}$ is defined to be zero if $v_i, v_j, v_k$ do not form a triangle). Then the area $A_{ijk}$ of $\Delta v_i v_j v_k$ satisfies the equation,

$$\theta_i^{jk} + \theta_j^{ik} + \theta_k^{ki} = \pi + \lambda A_{ijk}.$$

Since $a_i = \sum_{j,k} \theta_i^{jk}$, we have

$$\sum_i K_i = 2\pi N - \sum_{i,j,k} \theta_i^{jk}$$

$$= 2\pi N - \sum_{\Delta v_i v_j v_k} (\theta_i^{jk} + \theta_j^{ik} + \theta_k^{ki})$$

$$= 2\pi N - \pi|F| - \lambda Area(X)$$

$$= 2\pi\chi(X) - \lambda Area(X).$$

As a consequence, if $K^2 = E^2$, then the average curvature of a circle-packing metric is $2\pi\chi(X)/N$. In particular, the constant curvature metric has curvature $2\pi\chi(X)/N$.

3.2. *The Combinatorial Ricci Flow*

We define the combinatorial Ricci flow (with background metric $K^2$) to be

(3.1) $$dr_i(t)/dt = -K_i s(r_i(t)).$$

The evolution of curvature was first obtained by He [He] for zero weighted triangulations in Euclidean geometry.

**Proposition 3.2.** *Under the Ricci flow (3.1), the curvature $K_i(t)$ evolves according to*

(3.2) $$dK_i/dt = \sum_{j \sim i} C_{ij}(K_j - K_i) + \lambda B_i K_i$$

where the sum is over all vertices $v_j$ adjacent to $v_i$, $C_{ij}$ and $B_i$ are positive elementary functions in the radii $r_1, ..., r_N$. Furthermore, $C_{ij} = C_{ji}$.

This is a direct consequence of proposition 2.4 since $dK_i/dt = -\sum_{j,k} d\theta_i^{jk}/dt$ where $\theta_i^{jk}$ are inner angles at $v_i$ in the triangle $\Delta v_i v_j v_k$.

As a consequence, we obtain the maximum principle for curvature evolutions.

**Corollary 3.3.** (The Maximum Principle) *Let $r(t) = (r_1(t), ..., r_N(t))$ be a solution of the Ricci flow (3.1) in an interval. Let $M(t) = \max\{K_1(t), ..., K_N(t)\}$ and $m(t) = \min\{K_1(t), ..., K_N(t)\}$.*
*(1) In the case of Euclidean geometry $E^2$, then the function $M(t)$ is non-increasing in $t$ and $m(t)$ is non-decreasing in time $t$.*
*(2) In the case of hyperbolic geometry $H^2$, then $\max(M(t), 0)$ is non-increasing in $t$ and $\min(m(t), 0)$ is non-decreasing in time $t$.*

3.3. *The long time existence of solutions*

As a consequence of the maximum principle, we have,

**Proposition 3.4**. *Suppose the background geometry $K^2$ is $E^2$ or $H^2$. For any initial metric $r(0) \in \mathbf{R}_{>0}^N$, the solution to the Ricci flow (3.1) with the given initial metric exists for all time $t \geq 0$.*

Indeed, we note first of all the curvature $K_i$ is always bounded, i.e.,

$$(2-d)\pi < K_i < 2\pi$$

where $d$ is the degree of the vertex $v_i$.

In the case of Euclidean geometry $E^2$, by (3.1), we obtain

$$c_1 e^{(2-d)\pi t} \leq r_i(t) \leq c_2 e^{2\pi t}$$

for some positive constants $c_1$ and $c_2$. Thus $r_i(t)$ remains bounded above and away from zero as long as $t$ is bounded. This implies that the Ricci flow (3.1) has a solution for all time $t \geq 0$ for any initial value $r(0) \in \mathbf{R}_{>0}^N$.

In the case of $H^2$, since $K_i(t) = d(\ln(\coth(r_i/2)))/dt$, we obtain

$$\coth r_i/2 \leq c e^{2\pi t}$$

where $c$ is some constant. Thus $r_i(t)$ remains bounded away from 0 as long as time $t$ is bounded. On the other hand, the maximum principle (corollary 3.3) shows that $r_i(t)$ is bounded from above as well. To see this, we first have,

**Lemma 3.5.** *For any $\epsilon > 0$, there exists a number $L$ so that when $r_i > L$, the inner angle $\theta_i$ in the hyperbolic triangle $\Delta v_i v_j v_k$ is smaller than $\epsilon$.*

This is a simple fact from hyperbolic geometry. Indeed, if we put $v_i$ in the Euclidean center of the disk model for the hyperbolic plane and choose a large radius circle centered at $v_i$ in the hyperbolic space, then the Euclidean diameters of the other two circles centered at $v_j$ and $v_k$ can be made very small. Thus the angle at $v_i$ can be made arbitrary small if the radius $r_i$ is large.

As a consequence of the lemma, if there exists a sequence of time $t_n \to a$ where $a \in (0, \infty]$ so that $r_i(t_n)$ tends to infinity, then $K_i(t_n)$ will tend to $2\pi$. This contradicts the non-increasing property of the function $M(t) < 2\pi$.

By combining the above two, we see that as long as time remains bounded, then $r_i(t)$ and $1/r_i(t)$ remain bounded. This implies that the solution to (3.1) exists for all time $t \geq 0$.

As a consequence of lemma 3.5 and the maximum principle, we have,

**Corollary 3.6.** *Suppose $r(t)$ ($t \geq 0$) is a solution to the Ricci flow (3.1) in hyperbolic background geometry, then there exists a constant $C$ so that $r_i(t) \leq C$ for all time $t \geq 0$ and all indices $i$.*

3.4. *An existence theorem in the hyperbolic and Euclidean cases*

The following proposition shows that it suffices to establish a compactness result in order to show the long time convergence of the solution to (3.1).

**Proposition 3.7.** *Suppose $r(t)$ for $t \in [0, \infty)$ is a solution to the Ricci flow (3.1) in $H^2$ geometry so that the set $\{r(t) | t \in (0, \infty)\}$ lies in a compact region in $\mathbf{R}_{>0}^N$. Then $r(t)$ converges exponentially fast to a circle packing metric in $\mathbf{R}_{>0}^N$ whose curvature at each vertex is zero.*

**Proof.** Since the coefficient $B_i$ in (3.2) is an elementary function in $r_1, ..., r_N$ and is always positive, it follows from the assumption that there exist two positive constants $c_1$ and $c_2$ so that
$$c_1 \leq B_i(r_1(t), ..., r_N(t)) \leq c_2$$
for all time $t \geq 0$. Let $M_*(t) = \max(M(t), 0)$ and $m_*(t) = \min(m(t), 0)$. By equation (3.2), we obtain
$$dM_*(t)/dt \leq -c_1 M_*(t)$$
and
$$dm_*(t)/dt \geq -c_2 m_*(t).$$
Thus there are two constants $c_3 > 0$ and $c_4 > 0$ so that
$$c_5 e^{-c_3 t} \leq m_*(t) \leq M_*(t) \leq c_6 e^{-c_4 t}.$$

This shows that the curvature $K_i(t)$ converges exponentially fast to zero. This in turn implies that $\int_0^t K_i(s) ds$ converges exponentially fast to some constant. Since (3.1) says $\coth(r_i(t)/2) = ce^{\int_0^t K_i(s)ds}$, it follows that $\lim_{t\to\infty} r_i(t) = R_i$ exists in $\mathbf{R}_{>0}$ and $r_i(t)$ converges to $R_i$ exponentially fast. Thus the result follows. QED

The same proposition also works for the Euclidean case. In this case, we should normalize the equation. By the Gauss-Bonnet theorem, the sum of the curvatures is $2\pi\chi(X)$. Thus the average curvature $K_{av}$ is $2\pi\chi(X)/N$ which is independent of time $t$. We define the normalized Ricci flow to be

$$(3.3) \qquad dr_i(t)/dt = -(K_i - K_{av})r_i$$

It is easy to verify that $r(t)$ solves (3.3) for $K^2 = E^2$ if and only if $e^{K_{av}t}r(t)$ solves (3.3).

**Proposition 3.8.** *Suppose $r(t)$ for $t \in (0, \infty)$ is a solution to the normalized Ricci flow (3.3) in $E^2$ so that the set $\{r(t)|t \in (0,\infty)\}$ lies in a compact region in $\mathbf{R}_{>0}^N$. Then $r(t)$ converges exponentially fast to a circle packing metric in $\mathbf{R}_{>0}^N$ whose curvature at each vertex is the average curvature $K_{av}$.*

**Proof.** To show that the metrics $r(t) = (r_1(t), ..., r_N(t))$ converge exponentially to a vector in $\mathbf{R}_{>0}^N$, by the same argument as in the proof of proposition 3.7, it suffices to show that the curvatures $K(t) = (K_1(t), ..., K_N(t))$ converges exponentially fast to $K_{av}(1, ..., 1)$. To this end, let us consider the function $g(t) = \sum_{i=1}^N (K_i(t) - K_{av})^2$. Since $\sum_{i=1}^N K_i = NK_{av}$, we have $g(t) = \sum_{i=1}^N K_i^2 - NK_{av}^2$. Thus the derivative $g'(t) = 2\sum_{i=1}^N K_i(t)K_i'(t)$. By proposition 2.4 and $dK_{av}/dt = 0$, the evolution equation for the curvature $K_i(t)$ satisfies the following,

$$K_i'(t) = \sum_{j\sim i} C_{ij}((K_j - K_{av}) - (K_i - K_{av}))$$

$$= \sum_{j\sim i} C_{ij}(K_j - K_i)$$

where $C_{ij} = C_{ji}$ is positive and is an elementary function in $r_1, ..., r_N$. Thus

$$g'(t) = 2\sum_{i=1}^N \sum_{i\sim j} C_{ij}K_i(K_j - K_i).$$

By switching the order of $i, j$, we also have

$$g'(t) = 2\sum_{j=1}^N \sum_{j\sim i} C_{ji}K_j(K_i - K_j).$$

Thus

$$(3.4) \qquad g'(t) = -\sum_{i\sim j} C_{ij}(K_i - K_j)^2$$

where the sum is over all pairs of adjacent vertices $v_i$ and $v_j$. Since the metrics $r(t)$ lies in a compact region in $\mathbf{R}_{>0}^N$, the coefficients $C_{ij}(r_1(t), ..., r_N(t))$ remain bounded for all time $t \geq 0$, i.e., $C_{ij} \geq c_1 > 0$ for some constant $c_1$.

**Claim.** There exists a constant $c_2 > 0$ depending only on the triangulation $T$ so that

$$c_2 \sum_{i=1}^{N} (K_i - K_{av})^2 \leq \sum_{i \sim j} C_{ij} (K_i - K_j)^2$$

for all time $t \geq 0$.

Assuming the claim, by (3.4), we obtain the following

$$g'(t) \leq -c_2 g(t).$$

This implies that $g(t) \leq c_3 e^{-c_2 t}$ for all time $t \geq 0$ for some constant $c_3$. In particular, we obtain $|K_i(t) - K_{av}|^2 \leq g(t) \leq c_3 e^{-c_2 t}$. Thus the proposition holds.

To see the claim, we first note that since $C_{ij}$ are uniformly bounded from below, it suffices to prove that

(3.5)
$$\sum_{i=1}^{N} (K_i - K_{av})^2 \leq c_4 \sum_{i \sim j} (K_i - K_j)^2$$

for some constant $c_4$. By Cauchy's inequality that for any sequence $\{a_i\}$

$$(a_i - \frac{\sum_{j=1}^{N} a_j}{N})^2 \leq 1/N (\sum_{j=1}^{N} (a_i - a_j)^2),$$

we obtain

$$\sum_{i=1}^{N} (K_i - K_{av})^2 \leq 1/N \sum_{i,j=1}^{N} (K_i - K_j)^2$$

On the other hand, since the surface $X$ is connected, for any two vertices $v_i$ and $v_j$ in $T$, there exists a sequence of vertices $v_{m_1} = v_i, v_{m_2}, ..., v_{m_l} = v_j$ so that $v_{m_h}$ and $v_{m_{h+1}}$ are adjacent. We can estimate the value $(K_i - K_j)^2$ in terms of the sum of adjacent terms $(K_{m_h} - K_{m_{h+1}})^2$. This shows that there exists a constant $c_5$ depending only on $T$ so that $(K_i - K_j)^2 \leq c_5 \sum_{l \sim m} (K_l - K_m)^2$. Thus

$$\sum_{i,j=1}^{N} (K_i - K_j)^2 \leq c_5 N \sum_{i \sim j} (K_i - K_j)^2.$$

Combining all these, we establish the claim. Thus the proposition follows.

### 3.5. *The Ricci flow is variational*

In this section, we establish the fact that the combinatorial Ricci flow is essentially the gradient flow of a function. In the case that the background geometry is $E^2$ or $H^2$, the function is convex.

Let $u_i$ be a function of $r_i$ so that $du_i/dr_i = 1/s(r_i)$. To be more precise, we choose $u_i = \ln r_i$ when $K^2 = E^2$, $u_i = \ln \tanh(r_i/2)$ in $K^2 = H^2$ and $u_i = \ln \tan(r_i/2)$ when $K^2 = S^2$. Let $U$ be the following open set. If $K^2 = E^2$, then $U = \mathbf{R}^N$; if $K^2 = H^2$, then $U = (-\infty, 0)^N$. In the case of the spherical geometry, we take $U$ to be the inverse image of all possible $r = (r_1, ..., r_N)$ under the map $(u_1, ..., u_N) \to (r_1, ..., r_N)$.

Under this change of variable, the Ricci flow (3.1) takes the following form

$$(3.6) \qquad du_i/dt = -K_i.$$

Indeed, (3.6) is a gradient flow in the variable $u$. To see this, by lemma 2.3, we have $\partial K_i/\partial r_j s(r_j) = \partial K_j/\partial r_i s(r_i)$. This is equivalent to $\partial K_i/\partial u_j = \partial K_j/\partial u_i$. Thus the 1-form $\sum_{i=1}^N K_i du_i$ is closed in the simply connected space $U$. Define a function $f(u) = \int_a^u \sum_{i=1}^N K_i du_i$ where $a$ is any point in $U$. Then we have $\partial f(u)/\partial u_i = K_i$. This shows that (3.6) is the negative gradient flow of the function $f$. We remark that a family of functionals like this $f$ was first constructed by Colin de Verdiere in his paper [CV] in the case the weight function is zero. Colin de Verdiere also first proved the convexity of $f$ in the zero weight case.

**Proposition 3.9.** *(a) If $K^2 = H^2$, then the function $f(u) : U \to \mathbf{R}$ is strictly convex.*

*(b) If $K^2 = E^2$, then the function $f(u) : U \to \mathbf{R}$ is convex. Furthermore, it satisfies $f(u + a(1, 1, ..., 1)) = f(u)$ for all $a$ and $f$ is strictly convex when restricted in the plane $P = \{(u_1, ..., u_N) | \sum_{i=1}^N u_i = 0\}$.*

**Proof.** Let us consider the hessian of the function $f$.

Let $a_{ij} = \partial^2 f/\partial u_i \partial u_j = \partial K_i/\partial u_j$. By lemma 2.2, we have $a_{ii} > 0$ and $a_{ij} \leq 0$ ($a_{ij} < 0$ if $v_i$ is adjacent to $v_j$). Furthermore, lemmas 2.2 and 2.3 show that if $K^2 = H^2$, then $\partial \theta_i/\partial r_i \sinh(r_i) + \partial \theta_i/\partial r_j \sinh(r_j) + \partial \theta_i/\partial r_k \sinh(r_k) < 0$. Now sum over all such inequalities at the $i$-th vertex and use the fact that $a_{ij} = -\sum_{r,s} \partial \theta_i^{rs}/\partial r_j \sinh(r_j)$, we obtain $\sum_{j=1}^N a_{ij} > 0$. Recall a simple fact from linear algebra.

**Lemma 3.10.** *Suppose $A = [a_{ij}]_{n \times n}$ is a symmetric metrix.*

*(a) If $a_{ii} > \sum_{j \neq i} |a_{ij}|$ for all indices $i$, then $A$ is positive definite*

*(b) If $a_{ii} > 0$ and $a_{ij} \leq 0$ for all $i \neq j$ so that $\sum_{i=1}^n a_{ij} = 0$ for all $j$, then $A$ is semi-positive definite so that its kernel is one-dimensional.*

Thus part (a) of the proposition follows.

To see part (b) of the proposition, note that lemmas 2.2 and 2.3 show that $a_{ii} > 0$, $a_{ij} \leq 0$ and $\sum_{i=1}^N a_{ij} = 0$. By the linear algebra fact above, this implies that the matrix $[a_{ij}]$ is semi-definite. Thus $f$ is convex in $\mathbf{R}^n$. The identity $f(u + (a, a, ..., a)) = f(u)$ follows from the fact that $K_i(ar) = K_i(r)$. The last fact that $f|_P$ is strictly convex follows from the linear algebra consideration, namely, the kernel of the matrix $[a_{ij}]_{N\times N}$ is generated by the vector $(1,1,...1)$.

To prove the lemma 3.10 in the case of (a), it suffices to show that $\det(A) > 0$, since the same argument will show that all determinants of the principle submatrices of $A$ are positive. This in turn implies that $A$ is positive definite. Now to show $\det(A) > 0$, we first show that $\det(A) \neq 0$. This follows easily from the given condition. For instance, if $(x_1, ..., x_n)$ were a non-zero vector so that $\sum_j a_{ij} x_j = 0$ for all $i$, take $x_k$ to be the coordinate with the maximum absolute value. Now consider $\sum_{j=1}^n a_{kj} x_j = 0$. We obtain $a_{kk}|x_k| = |\sum_{i\neq k} a_{ki} x_i|$. But this contradicts the assumptions that $a_{kk} > \sum_{i\neq k} |a_{ki}|$ and $|x_k| \geq \max_i\{|x_i|\}$. Next multiply all off diagonal entries of $A$ by a constant $t$ so that $|t| \leq 1$ to get a new matrix $B(t)$. The matrix $B(t)$ satisfies the same conditions that $A$ satisfies, thus $\det(B(t)) \neq 0$. But $\det(B(0)) > 0$, thus $\det(B(t)) > 0$. In particular $\det(A) = \det(B(1)) > 0$.

To see part (b) of the lemma, first of all it is clear that $[a_{ij}]$ is semi-positive definite (since if we replace $a_{ii}$ by $a_{ii} + s$ for a small positive number $s$, it becomes positive definite by part (a)). Thus all eigenvalues of it are non-negative. Now it is easy to see that the kernel is one-dimensional spanned by $(1, 1, , 1)$. Thus the result follows. QED

Given a smooth strictly convex function defined on an open convex set $D$, the gradient map sending a point in $D$ to the gradient of $f$ at the point is always injective. Thus as a consequence of the above proposition, we obtained a new proof of the uniqueness part of Thurston's theorem.

**Corollary 3.11.** ([Th]) *Let* $\Pi : \mathbf{R}_{>0}^N \to \mathbf{R}^N$ *be the map sending a metric* $(r_1, ..., r_N)$ *to the corresponding curvature* $(K_1, ..., K_N)$ *where* $K^2 = H^2$ *or* $E^2$. *Then*

*(a) in the hyperbolic background geometry, $\Pi$ is injective, i.e., the metric is determined by its curvature.*

*(b) in the Euclidean geometry, $\Pi$ restricted to the subset $\{(r_1, ..., r_N) \in \mathbf{R}_{>0}^N | \Pi_{i=1}^N r_i = 1\}$ is injective, i.e., the metric is determined by its curvature up to a scalar multiplication.*

## §4. Degeneration of Circle Packing Metrics

We will recall a related work of [Th]. For completeness a slightly different proof is presented. The proof below follows closely [MR]. For simplicity, if $I = \{r_j, ...r_l\}$ is a subset of vertices $V$, we will also use $I$ to denote the index set $\{j, ..., l\}$.

**Proposition. 4.1** ([Th]) *Suppose $(T, \Phi)$ is a weighted generalized triangulation of a closed surface $X$ and $I$ is a subset of vertices $V$. Here the weight is a map $\Phi : E \to [0, \pi)$. Let $r^{(n)} = (r_1^{(n)}, ..., r_N^{(n)})$ be a sequence of circle packing metrics based on $(T, \Phi)$ in background geometry $K^2$ so that $\lim_{n\to\infty} r_i^{(n)} = 0$ for $i \in I$ and $\lim_{n\to\infty} r_i^{(n)} > 0$ for $i \notin I$. Then*

$$\lim_{n\to\infty} \sum_{i\in I} K_i(r^{(n)}) = - \sum_{(e,v)\in Lk(I)} (\pi - \Phi(e)) + 2\pi\chi(F_I),$$

where $F_I$ is the CW-subcomplex consisting of cells whose vertices are in $I$ and $Lk(I) = \{(e,v)|e$ is an edge so that $e \cap I = \emptyset$ and the vertex $v \in I$ so that $e,v$ form a triangle$\}$.

Furthermore, if the weight $\Phi : E \to [0, \pi/2]$ and the background geometry is $E^2$ or $H^2$, then for any circle packing metric $r$ based on $(T, \Phi)$ and any proper subset $I$ of vertices $V$, we have
$$\sum_{i \in I} K_i(r) > - \sum_{(e,v) \in Lk(I)} (\pi - \Phi(e)) + 2\pi\chi(F_I).$$

**Proof.** Recall that $\theta_i^{jk}$ denotes the inner angle at $v_i$ in the triangle $\Delta v_i v_j v_k$. Let $\tilde{X} \to X$ be the universal cover. Consider triangles in $T$ having a vertex in $I$. These triangles can be classified into three types $A_1$, $A_2$ and $A_3$ where a triangle $\Delta$ is in $A_i$ if and only if a lift of it in the universal cover has exactly $i$ many vertices in $\pi^{-1}(I)$. (Note that in the surface $X$, there may be triangles having all vertices being $v_i$.)

Now the curvature $K_i = 2\pi - a_i$ where $a_i$ is the cone angle at $v_i$. Here $\sum_{i \in I} a_i$ can be written as
$$\sum_{v_i \in I, \Delta v_i v_j v_k \in A_1} \theta_i^{jk} + \sum_{v_i, v_j \in I, \Delta v_i v_j v_k \in A_2} (\theta_i^{jk} + \theta_j^{ik}) + \sum_{\Delta v_i v_j v_k \in A_3} (\theta_i^{jk} + \theta_j^{ik} + \theta_k^{ij}).$$

It is known that as the radius $r_i^{(n)}$ tends to zero, the following holds.

(1) If $\Delta v_i v_j v_k \in A_3$, then $\lim_{n \to \infty}(\theta_i^{jk}(r^{(n)}) + \theta_j^{ik}(r^{(n)}) + \theta_k^{ij}(r^{(n)})) = \pi$. Indeed, the triangle is shrinking to a point in $K^2$ and thus is becoming a Euclidean triangle.

(2) If $\Delta v_i v_j v_k \in A_2$ and $i, j \in I$, then $\lim_{n \to \infty}(\theta_i^{jk}(r^{(n)}) + \theta_j^{ik}(r^{(n)})) = \pi$. Indeed, the triangle degenerates to a geodesic segment.

(3) If $\Delta v_i v_j v_k \in A_1$ and $i \in I$, then $\lim_{n \to \infty} \theta_i^{jk}(r^{(n)}) = \pi - \Phi(v_j v_k)$ by definition. Here $v_j v_k \in Lk(I)$.

To summary, we obtain

$$\lim_{n \to \infty} \sum_{i \in I} K_i(r^{(n)}) = - \sum_{(e,v) \in Lk(I)} (\pi - \Phi(e)) + 2\pi(|I| - 1/2|A_2| - 1/2|A_3|).$$

Here $|A_i|$ denotes the number of triangles of type $A_i$. Note $|A_3|$ is the number of 2-cells in $F_I$ and $|I|$ is the number of 0-cells in $F_I$. By the construction, we see that the number of edges in $F_I$ is $1/2(|A_2| + 3|A_3|)$. Thus $|I| - 1/2|A_2| - 1/2|A_3|$ is the Euler characteristic of $F_I$. This ends the proof.

To see the second part of the proposition, we note that since the background geometry is $H^2$ or $E^2$, the sum of the inner angles of a triangle is at most $\pi$. Thus, in the first limit (1) above, $\pi$ is the maximal value. By the same argument, in the second limit (2) above, $\pi$ is again the maximal value which is never achieved by a non-degenerate triangle. In the third limit (3), since $\Phi(e) \leq \pi/2$, using lemma 2.1 on the monotonicity, we see that $\pi - \Phi(v_j v_k)$ is also the maximal value for all possible circle packing metrics based on $(T, \Phi)$. Furthermore, since $I \neq V$, the cases (2) or (3) do occur. This shows that $-\sum_{(e,v) \in Lk(I)} (\pi - \Phi(e)) + 2\pi\chi(F_I)$ is strictly less than the sum of the curvature $\sum_{i \in I} K_i(r)$ for any circle packing metric based on $(T, \Phi)$. QED.

**Remark.** The proof works for all CW-cell decompositions so that each cell has 3 vertices and each vertex has degree at least 2.

**Corollary 4.2.** *Make the same assumption as in proposition 4.1 and assume that the weight $\Phi : E \to [0, \pi/2]$. Suppose $-\sum_{(e,v) \in Lk(I)} (\pi - \Phi(e)) + 2\pi\chi(F_I) > 0$. Then,*

*(a) either there exists a null homotopic loop consisting of 3 edges $e_1, e_2, e_3$ so that $\sum_{i=1}^{3} \Phi(e_i) > \pi$. Furthermore, if $X \neq S^2$, then these 3 edges do not form the boundary of a triangle, or*

*(b) $X = S^2$ and $I = V - \{v\}$ for one vertex $v$.*

Indeed, choose $I'$ be the non-empty subset of $I$ of the smallest cardinality so that $-\sum_{(e,v) \in Lk(I)} (\pi - \Phi(e)) + 2\pi\chi(F_I) > 0$. In this case, $F_{I'}$ is connected since the quantity is additive with respect to connected components of $F_I$. Thus $\chi(F_{I'}) \leq 1$ so that the equality holds if and only if $F_{I'}$ is contractible. Since $-\sum_{(e,v) \in Lk(I')} (\pi - \Phi(e)) \leq 0$, we see that $\chi(F_{I'}) = 1$ and $F_{I'}$ is contractible. Let $e_1, ..., e_k$ be the set of all edges in $Lk(I')$. They form a null homotopic loop due to the contractibility of $F_{I'}$. Now if $Lk(I') \neq \emptyset$, then $\sum_{(e,v) \in Lk(I')} (\pi - \Phi(e)) \geq k/2\pi$. Thus $k = 3$. Note $k = 1, 2$ are ruled out by the assumption on the generalized triangulation $T$. If these three edges form the boundary of a triangle, then $F_{I'}$ must be the complement of the triangle. Thus $X = S^2$.

If $Lk(I') = \emptyset$, then $V - I'$ must consists of vertices so that no two of these vertices are adjacent. On the other hand, $F_{I'}$ is contractible. Thus $V - I'$ consists of only one point. Thus we see that $X$ is a union of two contractible sets and $X = S^2$.

**Corollary 4.3.** *Make the same assumption as in proposition 4.1 and assume that the weight $\Phi : E \to [0, \pi/2]$. If there are no subsets of $V$ satisfying the condition in corollary 4.2 and if $-\sum_{(e,v) \in Lk(I)} (\pi - \Phi(e)) + 2\pi\chi(F_I) = 0$, then*

*(a) either there exists a null homotopic loop consisting of $k = 3$ or $4$ edges $e_1, ..., e_k$ so that $\sum_{i=1}^{k} \Phi(e_i) = (k - 2)\pi$. Furthermore, if $X \neq S^2$, then these $k$ edges do not form the boundary of a triangle, or union of two adjacent triangles, or*

*(b) $X = S^2$ and $I = V - \{v, v'\}$ for two non-adjacent vertices $v, v'$, or*

*(c) $X = \mathbf{R}P^2$ and $I = V - \{v\}$.*

Let $I'$ be the subset of the smallest cardinality so that $\sum_{(e,v) \in Lk(I)} (\pi - \Phi(e)) + 2\pi\chi(F_I) = 0$. By the assumption, we conclude that $F_{I'}$ is connected. Thus $\chi(F_{I'}) \leq 1$. By the hypothesis, $\chi(F_{I'}) \geq 0$. If $\chi(F_{I'}) = 1$, then by the given hypothesis, $F_{I'} \neq \emptyset$. By the same argument as above, we obtain conclusion (a). If $\chi(F_I) = 0$, then $Lk(I') = \emptyset$, i.e., $V - I'$ consists of vertices so that no two of them are adjacent. But $F_{I'}$ is either the annulus or the Möbuis band. This shows that either $X = S^2$ and $I' = V - \{v, v'\}$ where $v, v'$ are not adjacent or $X = \mathbf{R}P^2$ and $I' = V - \{v\}$.

**Corollary 4.4.** *Suppose $X$ is a closed surface of negative Euler characteristic and $(T, \Phi)$ is a generalized triangulation on $X$ so that the conditions (1.5a) and (1.5b) in theorem 1.2 hold. If $r^{(n)} \in \mathbf{R}_{>0}^N$ is a sequence of circle packing metrics in the hyperbolic background metric based on $(T, \Phi)$ so that $\liminf_{n \to \infty} K_i(r^{(n)}) \geq 0$ and $\limsup_{n \to \infty} K_i(r^{(n)}) < 2\pi$ for all indices $i$, then $r^{(n)}$ contains a convergent subsequence in $\mathbf{R}_{>0}^N$.*

Indeed, by lemma 3.5 and $\limsup_{n \to \infty} K_i(r^{(n)}) < 2\pi$, we see that $r_i^{(n)}$ is bounded from above. To see that $\{r_k^{(n)}\}$ is bounded away from 0, suppose otherwise that it contains a subsequence $\{r_k^{(n_j)}\}$ converging to zero. Let $I$ be the non-empty subset of indices so that $\lim_{j \to \infty} r_i^{(n_j)} = 0$ for $i \in I$ and $\lim_{j \to \infty} r_i^{(n_j)} > 0$ for $i \notin I$. Since $\liminf_{j \to \infty} K_i(r^{(n_j)}) \geq 0$ for all indices $i$, we obtain $\lim_{j \to \infty} \sum_{i \in I} K_i(r^{(n_j)}) \geq 0$. By corollaries 4.2 and 4.3, we see that one of the conditions (1.5a), (1.5b) in theorem 1.2 must be false. This contradicts the assumption.

### §5. A Proof of Theorem 1.2

We will prove the following result which is slightly more general than theorem 1.2.

**Theorem 5.1.** *For any weighted generalized triangulation $(T, \Phi)$ of a closed connected surface $X$ of negative Euler characteristic and any initial choice of circle-packing metric, the solution to (1.4) with the given initial value exists for all time and converges if and only if the following two conditions hold.*

*(1.5a) for any three edges $e_1, e_2, e_3$ forming a null homotopic loop in $X$, if $\sum_{i=1}^{3} \Phi(e_i) \geq \pi$, then these three edges form the boundary of a triangle of $T$, and*

*(1.5b) for any four edges $e_1, e_2, e_3, e_4$ forming a null homotopic loop in $X$, if $\sum_{i=1}^{4} \Phi(e_i) \geq 2\pi$, then these $e_i$'s form the boundary of the union of two adjacent triangles.*

*Furthermore, if it converges, then it converges exponentially fast to a hyperbolic metric on $X$ without any singularity.*

Note that if the weight $\Phi$ is zero, then conditions (1.5a) and (1.5b) hold automatically. Thus the above result implies the result of Andreev-Koebe-Thurston that for any triangulation of a closed surface of negative Euler characteristic, there exists a hyperbolic metric on the surface and a tangential circle packing in the metric so that the combinatorial pattern of tangency is given by the triangulation.

The proof of the theorem goes in two steps. In the first step, we give a new proof of Thurston's existence of circle packing metric using Ricci flow. This implies that the strictly convex function $f$ whose negative gradient is $du_i/dt = -K_i$ where $u_i = \ln \tanh(r_i/2)$ has a critical point. In particular, this implies that the function $f$ is bounded from below. On the other hand, corollary 4.4 shows that the function $f$ satisfies the Palais-Smale condition. Thus the negative gradient flow of $f$ converges to the critical point of $f$. Combining this with proposition 3.6, we see that all gradient lines converge exponentially fast to the critical point.

5.1. In this subsection, we produce a new proof of the existence of the zero curvature circle packing metric based on the weighted generalized triangulation $(T, \Phi)$.

To begin, let us choose the initial metric $r(0) = (r_1(0), ..., r_N(0))$ with all $r_i(0)$ sufficiently large so that it curvature at $v_i$ is close to $2\pi$. This is possible by lemma 3.5. In particular, we have $k_i(0) \geq 0$. Let $r(t) = (r_1(t), ..., r_N(t))$ $(t \geq 0)$ be a solution to the Ricci flow (3.1) in $H^2$ background metric with the initial value $r(0)$. Then by the maximum principle corollary 3.3, we have $k_i(t) \geq 0$ for all time $t \geq 0$. On the other hand, by the maximum principle again we have $\limsup_{t\to\infty} k_i(t) < 2\pi$. Thus by corollary 4.4, we see that the set $\{r(t)|t \in [0, \infty)\}$ lies in a compact region in $\mathbf{R}^N_{>0}$. Thus by proposition 3.7, the solution $r(t)$ converges exponentially fast to the zero curvature metric. In particular, this shows the zero curvature metric exists.

5.2. To show the convergence of the Ricci flow (3.1) in $H^2$ background metric in the general case, we use the result in subsection 3.5. Recall that if we consider the change of variable $u_i = \ln \tanh(r_i/2)$, then the Ricci flow becomes $du_i/dt = -K_i$ which is proved to be a negative gradient flow. Let $f$ be a potential function. It is known by proposition 3.9 that the function $f$ is strictly convex defined on $(-\infty, 0)^N$. We claim that $f$ satisfies the Palais-Smale condition. As a consequence, we see all gradient lines converge to the critical point. Thus the result follows.

To verify the Palais-Smale condition, first of all by §5.1, we see that $f$ has a critical point which must be the minimal point due to the convexity. Thus $f$ is bounded from below. To verify the second condition in Palais-Smale, take a sequence of points $\{u^{(n)}\}$ in $(-\infty, 0)^N$ so that the gradient of $f$ at $u^{(n)}$ converges to zero. Let the corresponding points of $u^{(n)}$ be $r^{(n)} \in \mathbf{R}^N_{>0}$. Then we have that the curvature $K_i(r^{(n)})$ tends to zero for all indices $i$ as $n$ tends to infinity. By corollary 4.4, this implies that the set $\{r^{(n)}|n \in \mathbf{Z}_{>0}\}$ lies in a compact set. This in turn implies that $\{u^{(n)}\}$ contains a convergent subsequence in $(-\infty, 0)^N$. This ends the proof.

Finally the necessity of the conditions (1.5a) and (1.5b) which was proved in [Th] follows from proposition 4.1. Indeed, since the zero curvature metric exists, we conclude that for any subset $I \in V$, the quantity $-\sum_{(e,v)\in Lk(I)}(\pi - \Phi(e)) + 2\pi\chi(F_I) < 0$. On the other hand, if (1.5a) or (1.5b) is false, then take $I$ to be the vertices in the disk bounded by the null homotopic loop formed by the edges. We see that $-\sum_{(e,v)\in Lk(I)}(\pi - \Phi(e)) + 2\pi(\chi(F_I)) \geq 0$.

**Remark 5.3.** The proof in fact shows the following for a generalization of the Ricci flow in hyperbolic background metric. Namely, suppose there exists a circle packing metric $R$ so that it curvature is $C_i$ at the i-th vertex. Then the following flow $dr_i/dt = -(K_i - C_i)\sinh(r_i)$ has a solution for all time $t \geq 0$ for all initial metrics and the solution converges exponentially fast to $R$.

§5. **A Proof of Theorem 1.1 and Proposition 1.3**

We will prove the following theorem for generalized triangulations which implies theorem 1.1.

**Theorem 6.1.** *For any weighted generalized triangulation $T$ of a closed connected surface $X$ and for any initial circle-packing metric based on the weighted triangulation, the solution to the normalized combinatorial Ricci flow (1.2) in Euclidean background metric with the given initial value converges if and only if for any proper subset $I \subset V$,*

$$(1.3) \qquad 2\pi\chi(X)|I|/N > -\sum_{(e,v)\in Lk(I)} (\pi - \Phi(e)) + 2\pi\chi(F_I).$$

*Furthermore, if the solution converges, it converges exponentially fast to the metric of constant curvature $2\pi\chi(X)/N$.*

The proof is essentially the same as in §5 and goes in two steps. In the first step, we show that the constant curvature circle packing metric exists due to condition (1.3). Next we show that the normalized Ricci flow is the negative gradient flow of a convex function in the variable $u_i = \ln r_i$ and that the function $f$ satisfies the Palais-Smale condition. Combining this with proposition 3.6, we see that all gradient lines converge exponentially fast to the critical point.

6.1. In this subsection, we produce a proof of the existence of the constant curvature circle packing metric based on the weighted generalized triangulation $(T, \Phi)$ under condition (1.3). This result is implicitly in [Th]. Our proof follows the work [MR]. Let $P$ be the set $\{r = (r_1, ..., r_N) \in \mathbf{R}^N_{>0} | \prod_{i=1}^N r_i = 1\}$. Consider the curvature map $\Pi : P \to \mathbf{R}^N$ sending the metric $r$ to its curvature $\Pi(r) = (K_1(r), ..., K_N(r))$. By corollary 3.9, it is injective. Furthermore, by Gauss-Bonnet, its image lies in the hyperplane defined by the equation $\sum_{i=1}^N K_i = 2\pi\chi(X)$. Consider the convex polytope $Y$ in the hyperplane defined by the set of all linear inequalities $\sum_{i\in I} K_i > -\sum_{i\in I}(\pi - \Phi(e)) + 2\pi\chi(F_I)$ where $I$ is a proper subset of vertices. Then by proposition 4.1 the map $\Pi : P \to Y$ sends compact sets to compact sets. Since both $P$ and $Y$ are homeomorphic to $\mathbf{R}^{N-1}$, by the invariance of domain theorem, the map $\Pi : P \to Y$ is a homeomorphism. In particular, we see that the constant curvature metric exists if and only if (1.3) holds for all proper subsets $I$ of $V$.

6.2. To show the convergence of the Ricci flow for all initial choice of metrics, let us take $u_i = \ln r_i$. Then the normalized Ricci flow becomes $du_i/dt = -(K_i - K_{av})$. This later equation is the negative gradient flow of a convex function $f(u)$ defined on the u-space $\mathbf{R}^N$. Furthermore, if we restrict to the subspace $Z = \{u = (u_1, ..., u_N) \in \mathbf{R}^N | \sum_{i=1}^N u_i = 0\}$, the function $f$ is strictly convex proposition 3.8. By the result in subsection 6.1, we see that $f$ has a critical point. Thus it is bounded from below. We claim that $f$ satisfies the Palais-Smale condition on the plane $Z$. To this end, take a sequence of vectors $u^{(n)}$ in $Z$ so that the gradient of $f$ at $u^{(n)}$ tends to zero. Let $r^{(n)}$ be the corresponding points in the space $P$. Then $\lim_{n\to\infty} K_i(r^{(n)}) = K_{av}$ for all indices $i$. By the assumption that $K_{av}(1, ..., 1)$ lies in the open cell $\Pi(P)$, we see that $r^{(n)}$ lies in a compact set in $P$. Thus $\{u^{(n)}\}$ contains a convergent subsequence. This shows all gradient lines converge to the constant curvature metric. Thus by proposition 3.6, the solution $r(t)$ converges exponentially fast to the constant curvature metric.

6.3. The proof in fact shows the following for a generalization of the Ricci flow in Euclidean background metric. Namely, suppose there exists a circle packing metric $R$ so that it curvature is $C_i$ at the i-th vertex. Then the following flow $dr_i/dt = -(K_i - C_i)r_i$ has a solution for all time $t \geq 0$ for all initial metrics and the solution converges exponentially fast to $R$.

6.4. *A proof of proposition 1.3*

**Proposition 1.3** *(a) If the Euler characteristic of the surface $X$ is positive, then condition (1.3) is valid for all generalized triangulations of $X$ so that there are no null homotopic loops consisting of 3 edges $e_1, ..., e_3$ so that $\sum_{i=1}^{k} \Phi(e_i) \geq \pi$.*
*(b) If the Euler characteristic of the surface $X$ is zero, then (1.3) is valid for all triangulations of $X$ so that (1.5a) and (1.5b) hold.*

This follows from corollaries 4.2 and 4.3.

### §7. Some Questions

There are several questions which we find interesting relating to the results in the paper.

1. *Generalizing Lemma 2.2*
It can be shown that lemma 2.2 is false if the assumption $\Phi_i \in [0, \pi/2]$ is relaxed to $\Phi_i$'s in $[0, \pi)$. In fact there exists a choice of $\Phi_n \in [0, \pi)$ so that $\partial\theta_i/\partial r_i < 0$. But on the other hand, it seems for applications, the convexity of the function $f$ in proposition 3.9 is essential. Thus the question is for what choices of $\Phi_n$'s in $[0, \pi)$ is the matrix $[s(r_j)\partial\theta_i/\partial r_j]_{3\times 3}$ positive or semi-positive definite.

2. *2-Dimensional Combinatorial Ricci flow for all piecewise constant curvature metrics*
Hamilton's ricci flow works for all Riemannian metrics. Fix a triangulation on a surface, the analog of Riemannian metric seems to be the piecewise constant curvature metric obtained by assigning lengths to each edge in the triangulation. Is there any Ricci flow defined on the space of all these metrics? This may have some applications to the problems like what is the space of all geometric triangulations of a triangle with the fixed combinatorial data.

3. Investigate the 2-dimensional combinatorial Ricci flow in theorems 1.1 and 1.2 so that the combinatorial conditions (1.3) and (1.5a) and (1.6b) are not valid. The solution still exists for all time $t \geq 0$. Is it true that for any solution $r(t) = (r_1(t), ..., r_N(t))$ where $t \geq 0$, the limit $\lim_{t\to\infty} r_i(t)$ always exists in the extended set $[0, \infty]$?
In the case of the Euclidean background metric, for negative Euler characteristic closed surfaces, if $r(t)$ with $t \geq 0$ is a solution to the Ricci flow (1.3), is it true that the metric $r(t)$ will have non-positive curvature when $t$ is large?

4. Investigate the 2-dimensional combinatorial Ricci flow in the case that the weight is in $[0, \pi)$. For instance, can one produce a new proof of the results of Igor Rivin [Ri] using heat equations?

5. We showed that rigidity of the circle packing follows from the strictly convexity of the potential function in this paper (corollary 3.11). The convexity is based on the main technical lemma 2.3. Derive the main result of He [He] on the rigidity of locally finite infinite circle packings from the convexity of the potential function.

6. Is there any 3-dimensional combinatorial analogous of Hamilton's 3-dimensional Ricci flow [Ha2]?

The recent work of David Glickenstein's UCSD thesis work ([Gl]) seems to indicate that the analog of the Yamabe flow exists in the combinatorial setting in dimension 3 (see also [CR]).

### Appendix. Some computations involving triangles

*A1. Basic results on triangles*

Given a triangle $\Delta v_i v_j v_k$ in one of the three geometries $K^2 = S^2, E^2$ and $H^2$, let the inner angle at the $v_i$ vertex be $\theta_i$ and the length of the edge $v_j v_k$ be $x_i$. Recall that the function $s(x)$ is defined according to the geometry $K^2$ as follows: $s(x) = x$ for $E^2$, $s(x) = \sin(x)$ for $S^2$ and $s(x) = \sinh(x)$ for $H^2$. Let $A_{ijk} = s(x_i)s(x_j)\sin(\theta_k)$. By the sine law, the quantity $A_{ijk}$ is symmetric in $i, j, k$. The following is a very useful lemma for computing derivatives.

**Lemma A-1.** *In all three geometries we have,*
(a) $\partial \theta_i / \partial x_i = s(x_i)/A_{ijk}$.
(b) $\partial \theta_i / \partial x_j = -\partial \theta_i / \partial x_i \cos(\theta_k)$.
(c) $\partial \theta_i / \partial x_j s(x_i)$ *is symmetric in* $i, j$.
(d) *Consider* $x_k$ *as a function of* $x_i, x_j, \theta_i$. *Then*

$$\partial x_k / \partial x_i = \cos(\theta_j).$$

The proof is a simple application of the cosine law. We will omit the details. It is interesting to note that these derivatives formulas are the same for all three geometries.

One interesting consequence of lemma A-1(c) is the following,

**Corollary A-2.** *For any triangle in the geometry* $K^2$, *the 1-form* $\sum_{l \in \{i,j,k\}} \theta_l s(x_l) dx_l$ *in variables* $x_i, x_j, x_k$ *is closed.*

*A2. A proof of lemma 1.3*

A simpler proof using mathematical software Maple VII is given in the next subsection. In this subsection, we present a conventional proof.

Suppose we have a collection of three closed disks $C_i, C_j$, and $C_k$ centered at $v_i, v_j$ and $v_k$ so that their radii are $r_i, r_j$ and $r_k$ and their exterior intersection angles are $\phi_i, \phi_j, \phi_k \in [0, \pi)$ in geometry $K^2$. Consider the triangle formed by $v_i, v_j, v_k$. Let its inner angle at $v_n$ be $\theta_n$ and the edge length opposite to $v_n$ be $x_n$. We want to show that $\partial \theta_i / \partial r_j s(r_j)$ is symmetric in $i, j$.

First of all if $L_i$, $L_j$ and $L_k$ are the geodesic lines passing through the pairs of the intersection points of $\{C_k, C_j\}$, $\{C_k, C_i\}$ and $\{C_i, C_j\}$ ($L_i$ is the tangent to $C_j$ and $C_k$ if they intersect in one point), then the three lines $L_i, L_j, L_k$ intersect in one point.

To see this, we note that the result is well known in the Euclidean geometry $E^2$. To see the result in hyperbolic geometry, let $O$ be the point of intersection of $L_i$ and $L_j$. Take the disk model for $H^2$ so that $O$ is the Euclidean center. Then both $L_i$ and $L_j$ are Euclidean line segments. Let $L'_k$ be the Euclidean diameter passing through $O$ and one of the intersection point of $C_i$ with $C_j$. The by the result in the Euclidean case, we see that $L_k$ must contain the other point of intersection of $C_i$ and $C_j$. This shows, $L'_k = L_k$. The same argument also works for the spherical case where we consider the stereographic projection from the antipodal point of $O$. Let $u_r$ be an intersection point of the circles $C_s$ and $C_t$ where $\{r, s, t\} = \{i, j, k\}$.

To begin the proof, by lemma A1, we have

$$s(r_j)\partial\theta_i/\partial r_j = s(r_j)(\partial\theta_i/\partial x_i \partial x_i/\partial r_j + \partial\theta_i/\partial x_k \partial x_k/\partial r_j)$$

$$= s(r_j)\partial\theta_i/\partial x_i(\partial x_i/\partial r_j - \cos(\theta_j)\partial x_k/\partial r_j)$$

$$= \frac{s(x_i)\sin(\theta_j)}{A_{ijk}} \frac{s(r_j)(\cos(\theta_j^{ki'}) - \cos(\theta_j)\cos(\theta_j^{ik'}))}{sin(\theta_j)}$$

Here $\theta_j^{ik'}$ and $\theta_j^{ki'}$ are the inner angles at $v_j$ in the triangles $\Delta v_j v_i u_k$ and $\Delta v_j v_k u_i$. Note that the first fraction $s(x_i)\sin(\theta_j)/A_{ijk}$ is symmetric in $i, j$ due to the sine law. Thus it suffices to show that
$$\frac{s(r_j)(\cos(\theta_j^{ki'}) - \cos(\theta_j)\cos(\theta_j^{ik'}))}{sin(\theta_j)}$$

is symmetric in $i, j$.

To this end, we will need the following lemma about quadrilaterals in $K^2$ so that two opposite angles are right angles.

**Lemma A2.** *Suppose $Q$ is an immersed quadrilateral in geometry $K^2$ so that its four edges are $T, X, Y, Z$ counted clockwise. Suppose the inner angle between sides $\{T, X\}$ and $\{Y, Z\}$ are right angles and the inner angle between sides $\{X, Y\}$ is $\theta$, and edges $T, Y$ are disjoint. Let the lengths of sides $X, Y, Z$ be $x, y, z$ and $t$ be the length of $T$ of $X \cap Z = \emptyset$ and $t$ be the negative of the length of $T$ if $X$ intersects $Z$. Then the following holds.*

(a). In the Euclidean geometry $E^2$,

$$t = (y - x\cos\theta)/\sin\theta.$$

(b). In the hyperbolic geometry $H^2$,

$$\tanh(t) = \cosh(x)(\tanh(y) - \tanh(x)\cos(\theta))/\sin(\theta).$$

(c). In the spherical geometry $S^2$,

$$\tan(t) = \cos(x)(\tan(y) - \tan(x)\cos(\theta))/\sin(\theta).$$

The proof is a direct application of trigonometries for right angled triangles and quadrilaterals with three right angles. (See for instance [Vi]).

Another simple fact we shall use is the following formulas for right angled triangles. Suppose a right angled triangle has three edge lengths $a, b, c$ and whose opposite angles are $\alpha, \beta$ and $\pi/2$. Then we have, $\cos(\alpha) = b/c$ in $E^2$, $\cos(\alpha) = \tanh(b)/\tanh(c)$ in $H^2$, and $\cos(\alpha) = \tan(b)/\tan(c)$ in $S^2$.

Let the distance from the common intersection $O$ to the intersection point of the lines $L_k$ and $v_i v_j$ be $d$. Let $D$ be $d$ if $O$ is inside the triangle and $-d$ otherwise. It is symmetric in $i, j$. Then by the above lemma and the trigonometric formulas, we obtain that the expression $s(r_j)(\cos(\theta_j^{ki'}) - \cos(\theta_j)\cos(\theta_j^{ik'}))/sin(\theta_j)$ is equal to the $D$, $\tanh(D)$ and $\tan(D)$ respectively. This verifies the symmetry. Note that if all angles $\Phi_i$'s are in $[0, \pi/2]$, then the common intersection point $O$ is inside the triangle. Thus $d > 0$ which is what lemma 2.2(b) states.

### A3. *A proof of lemma 2.3 using Maple VII*

The simplest proof of the lemma is to use some mathematical software so that we do not have to go over all possible configurations. We will present such a proof using Maple VII below.

Let $T_3$ be the space of all possible configurations of pairwise intersecting three circles in $K^2$ up to isometry so that their centers form a triangle. The we claim that $T_3$ is connected. To see this, since the space of all triangles up to isometry is connected, it suffices to show that the subspace of all pairwise intersecting three circles with fixed centers is connected. To see the later space, let us exam the three radii. We can choose the radius of the first circle arbitrary as long as it lies in its natural interval (e.g, it is $(0, \infty)$ for $E^2$ and $H^2$, and is $(0, \pi/2)$ for $S^2$). Now the radius of the second circle must also be in an interval since it intersects the first circle. The radius of the third must be in an interval as well since it intersects the first and the second circle and these two circles intersect. This shows that the space is connected.

Now to show $\partial \theta_i / \partial r_j s(r_j)$ is symmetric in $i, j$, it is much convenient for the machine to verify that $(\partial \theta_i / \partial r_j s(r_j))^2$ is symmetric in $i, j$. On the other hand, both functions $\partial \theta_i / \partial r_j s(r_j)$ and $\partial \theta_j / \partial r_i s(r_i)$ are analytic in variables $\Phi_i, \Phi_j, \Phi_k, r_i, r_j, r_k$ and defined on a connected space so that they coincide in the subspace $r_i = r_j$. Thus $\partial \theta_i / \partial r_j s(r_j)$ must be symmetric in $i, j$ if its square is symmetric.

To verify the symmetry of $((\partial \theta_i / \partial r_j) r_j)^2$ in $i, j$ for Euclidean geometry $E^2$ using Maple VII, we run the following,

```
> yi := sqrt(rj^2 + rk^2 + 2*rj*rk*bi);
> yj := sqrt(ri^2 + rk^2 + 2*ri*rk*bj);
> yk := sqrt(rj^2 + ri^2 + 2*rj*ri*bk);
> ai := arccos(((yj)^2 + (yk)^2 - (yi)^2)/(2*(yj)*(yk)));
> fij := diff(ai, rj)*rj;
> g1 := simplify(fij);
> g2 := g1^2;
> g3 := subs({ri = rj, rj = ri, bi = bj, bj = bi}, g2) - g2;
> simplify(g3);    0
```

To verify the symmetry of $((\partial \theta_i/\partial r_j)\sinh(r_j))^2$ in $i,j$ for hyperbolic geometry $H^2$ using Maple VII, we run the following,

```
> yi := arccosh(cosh(rj)*cosh(rk) + sinh(rj)*sinh(rk)*bi);
> yj := arccosh(cosh(rk)*cosh(ri) + sinh(ri)*sinh(rk)*bj);
> yk := arccosh(cosh(rj)*cosh(ri) + sinh(rj)*sinh(ri)*bk);
> ai := arccos((cosh(yj)*cosh(yk) - cosh(yi))/(sinh(yj)*sinh(yk)));
> fij := diff(ai, rj)*sinh(rj);
> g1 := simplify(fij);
> g2 := g1^2;
> g3 := subs({ri = rj, rj = ri, bi = bj, bj = bi}, g2) - g2;
> simplify(g3);    0
```

To verify the symmetry of $((\partial \theta_i/\partial r_j)r_j)^2$ in $i,j$ for spherical geometry $S^2$ using Maple VII, we run the following,

```
> yi := arccos(cos(rj)*cos(rk) + sin(rj)*sin(rk)*bi);
> yj := arccos(cos(rk)*cos(ri) + sin(ri)*sin(rk)*bj);
> yk := arccos(cos(rj)*cos(ri) + sin(rj)*sin(ri)*bk);
> ai := arccos((-cos(yj)*cos(yk) + cos(yi))/(sin(yj)*sin(yk)));
> fij := diff(ai, rj)*sin(rj);
> g1 := simplify(fij);
> g2 := g1^2;
> g3 := subs({ri = rj, rj = ri, bi = bj, bj = bi}, g2) - g2;
> simplify(g3);    0
```

This ends the proof.

Bennett Chow
Dept. of Math., UCSD, La Jolla, CA 92093
email: benchow@math.ucsd.edu

Feng Luo
Dept. of Math., Rutgers Univ., New Brunswick, NJ 08854
email: fluo@math.rutgers.edu